\documentclass{amsart}

\usepackage{amsmath,amssymb,amsfonts,amstext,amsthm,bm}
\usepackage{url}
\usepackage{graphicx}
\usepackage{hyperref}
\usepackage{tikz}
\usepackage[all,cmtip]{xy} 
\usepackage{multicol}
\usepackage{verbatim}
\input xy
\xyoption{all}
\usepackage[all]{xy}
\usepackage{xcolor}
\usepackage{enumitem}

\newcommand{\nwc}{\newcommand}

\nwc{\aaa}{\mathcal{F}}
\nwc{\aap}{\mathcal{F}_{P}}
\nwc{\al}{\alpha}

\nwc{\C}{\mathbb{C}}
\nwc{\cb}{\overline{C}}
\nwc{\ccc}{\mathcal{C}}
\nwc{\ch}{\widehat{C}}
\nwc{\cin}{\textbf{(v)}}
\nwc{\cl}{C'}
\nwc{\cp}{\mathcal{C}_{P}}
\nwc{\cpll}{\mathfrak{c}_{P'}}
\nwc{\ct}{\widetilde{C}}

\nwc{\dd}{\mathcal{L}}
\nwc{\ddd}{\mathfrak{d}}
\nwc{\ddl}{\mathcal{L}'}
\nwc{\dlp}{\delta_{P}}
\nwc{\doi}{\textbf{(ii)}}

\nwc{\enq}{$$}

\nwc{\fl}{\flushleft}
\nwc{\fff}{\mathcal{F}}
\nwc{\ffp}{\mathcal{F}_{P}}
\nwc{\ffq}{\mathcal{F}_{Q}}
\nwc{\ffl}{\mathcal{F}'}

\nwc{\G}{\mathcal{G}}
\nwc{\Ga}{\Gamma}
\nwc{\gon}{{\rm gon}}
\nwc{\gtl}{\widetilde{g}}

\nwc{\hra}{\hookrightarrow}
\nwc{\hua}{h^{1}(C,\aaa )}

\nwc{\kk}{{\rm K}}
\nwc{\kp}{{\kappa}}

\nwc{\llb}{\mathcal{L}}

\nwc{\mb}{\mathbb}
\nwc{\mc}{\mathcal}
\nwc{\mm}{\mathfrak{m}}
\nwc{\mmp}{\mathfrak{m}_{P}}
\nwc{\mpd}{\mathfrak{m}_{P}^{2}}

\nwc{\nn}{\mathbb{N}}

\nwc{\ob}{\overline{\mathcal{O}}}
\nwc{\obr}{\mathcal{O}^*}
\nwc{\obp}{\overline{\mathcal{O}}_P}
\nwc{\och}{\mathcal{O}_{\hat{C}}}
\nwc{\oh}{\widehat{\mathcal{O}}}
\nwc{\ohp}{\widehat{\mathcal{O}}_{P}}
\nwc{\ol}{\mathcal{O}'}
\nwc{\oma}{\Omega (\mathfrak{a})}
\nwc{\omo}{\Omega (\mathcal{O})}
\nwc{\oo}{\mathcal{O}}
\nwc{\op}{\mathcal{O}_P}
\nwc{\opc}{\mathcal{O}_{P,C}}
\nwc{\oph}{\hat{\mathcal{O}}_{P}}
\nwc{\opl}{\mathcal{O}_{P}'}
\nwc{\oplc}{\mathcal{O}_{P,C}'}
\nwc{\opll}{\mathcal{O}_{P'}}
\nwc{\opt}{\tilde{\mathcal{O}}_{P}}
\nwc{\optt}{{\mathcal{O}}_{\tilde{P}}}
\nwc{\oq}{\mathcal{O}_{Q}}
\nwc{\oqt}{\tilde{\mathcal{O}}_{Q}}
\nwc{\ot}{\widetilde{\mathcal{O}}}
\nwc{\overop}{\bar{\oo}_{P}}

\nwc{\pb}{\overline{P}}
\nwc{\pbb}{P^*}
\nwc{\pbi}{\overline{P_{i}}}
\nwc{\pbr}{\overline{P_{r}}}
\nwc{\pgmd}{\mathbb{P}^{g+2}}
\nwc{\pgmu}{\mathbb{P}^{g+1}}
\nwc{\ph}{\hat{P}}
\nwc{\pp}{\mathbb{P}}
\nwc{\ppn}{\mathbb{P}^{n}}
\nwc{\prv}{\noindent\textbook{Proof}:}
\nwc{\pt}{\widetilde{P}}
\nwc{\ptl}{\tilde{P}}
\nwc{\pum}{\mathbb{P}^{1}}

\nwc{\qh}{\hat{Q}}
\nwc{\qtl}{\tilde{Q}}
\nwc{\qua}{\textbf{(iv)}}

\nwc{\ra}{\rightarrow}
\nwc{\rh}{\hat{R}}

\nwc{\sei}{\textbf{(vi)}}
\nwc{\sep}{\beq\ast\ \ast\ \ast\enq}
\nwc{\sig}{\sigma}
\nwc{\Sig}{\Sigma}
\nwc{\ssp}{{\rm S}_{P}}
\nwc{\sss}{{\rm S}}
\nwc{\sssh}{\widehat{{\rm S}}}
\nwc{\sys}{\mathcal{L}}

\nwc{\tre}{\textbf{(iii)}}

\nwc{\um}{\textbf{(i)}}

\nwc{\val}{\mathcal{V}}
\nwc{\vpb}{v_{\overline{P}}}
\nwc{\vtxp}{\widetilde{V}_{x,P}}
\nwc{\vv}{\mathcal{W}}
\nwc{\vvp}{\mathcal{W}_{P}}
\nwc{\vxp}{V_{x,P}}

\nwc{\wh}{\hat{\omega}}
\nwc{\whp}{\hat{\omega}_{P}}
\nwc{\woch}{\omega\cdot\mathcal{O}_{\hat{C}}}
\nwc{\woh}{\omega\cdot\hat{\mathcal{O}}}
\nwc{\ww}{\omega}
\nwc{\wwb}{\omega^*}
\nwc{\wwct}{\omega _{\widetilde{C}}}
\nwc{\wwh}{\widehat{\omega}}
\nwc{\wwhp}{\widehat{\omega}_P}
\nwc{\wwp}{\omega _{P}}
\nwc{\wwt}{\widetilde{\omega}}
\nwc{\wwtp}{\widetilde{\omega}_P}

\nwc{\zz}{\mathbb{Z}}

\def\pp{{\mathbb P}}
\def\ww{\omega}
\def\sys{\mathcal{L}}
\def\kao{\overline{\kappa}}
\def\oo{\mathcal{O}}
\def\cb{\overline{C}}
\def\fff{\mathcal{F}}
\def\ccc{\mathcal{C}}
\def\vv{\mathcal{W}}
\def\oph{\widehat{\oo}_P}
\def\op{\mathcal{O}_P}
\def\obp{\overline{\mathcal{O}}_P}
\def\ch{\widehat{C}}
\def\aaa{\mathcal{F}}
\def\dd{\mathcal{L}}

\newtheorem{coro}{Corollary}[section]

\newtheorem{thm}[coro]{Theorem}

\newtheorem{conj}[coro]{Conjecture}
\newtheorem{spec}[coro]{Speculation}
\newtheorem{ex}[coro]{Example}
\let \fl=\flushleft

\let \sub=\subset

\let \al=\alpha
\let \pr=\prime

\let \sub=\subset
\let \mb=\mathbb
\let \mc=\mathcal

\let \ra=\rightarrow

\usepackage{faktor}

\title{Towards Brill--Noether theory for cuspidal curves}
\date{\today}
\author{Ethan Cotterill}

\address{Instituto de Matem\'atica, Estat\'istica, e Computa\c{c}\~ao Cient\'ifica, UNICAMP, Rua S\'ergio Buarque de Holanda, 651, 13.083-859 Campinas SP, Brazil}

\email{cotterill.ethan@gmail.com}

\author{Renato Vidal Martins}
\address{Departamento de Matem\'atica, ICEx, UFMG,
Av. Ant\^onio Carlos 6627,
30123-970 Belo Horizonte MG, Brazil}
\email{renato@mat.ufmg.br}

\begin{document}

\maketitle 

\section{Abstract}
Understanding when an abstract complex curve of given genus comes equipped with a map of fixed degree to a projective space of fixed dimension is a foundational question; and Brill–Noether theory addresses this question via linear series, which algebraically codify maps to projective targets. Classical Brill--Noether theory, which focuses on smooth curves, has been intensively explored; but much less is known for singular curves, particularly for those with non-nodal singularities. In a one-parameter family of smooth curves specializing to a singular curve $C_0$, one expects certain aspects of the global geometry of the smooth fibers to ``specialize" to the local geometry of the singularities of $C_0$. Making this expectation quantitatively precise involves analyzing the arithmetic and combinatorics of semigroups ${\rm S}$ attached to discrete valuations defined on (the local rings of) these singularities. In this largely-expository note we focus primarily on Brill--Noether-type results for curves with {\it cusps}, i.e., unibranch singularities; in this setting, the associated semigroups are {\it numerical} semigroups with finite complement in $\mb{N}$.

\section{Introduction}
We now introduce the main players in the sequel. Throughout we work over $\mb{C}$, though all results should also hold over algebraically closed fields of sufficiently large positive characteristic. 

\subsection{Cuspidal rational curves} {\it Rational curves}, by which we mean projective algebraic curves of geometric genus zero, are among the simplest algebraic varieties; yet via their singularities, their geometry intertwines with quantum topology and representation theory.

\medskip
It is often convenient to identify a rational curve of degree $d$ in $\mb{P}^n$ with its normalization, given by an $(n+1)$-tuple of homogeneous polynomials of degree $d$ over $\mb{P}^1$. Formally inverting one of these in a point where it is nonzero yields an $n$-tuple of power series $f_i(t)$, $i=1,\dots,n$ in a uniformizing parameter $t$; and any unibranch singularity, or {\it cusp}, is locally modeled in this way. Equivalently, its local algebra is the image of a ring map
\begin{equation}\label{ring_map}
f: R= \mb{C}[x_1,\dots,x_n] \ra \mb{C}[\![t]\!]
\end{equation}
defined by $x_i \mapsto f_i(t)$, $i=1,\dots,n$. Letting $v_t$ denote the standard $t$-adic valuation on $\mb{C}[\![t]\!]$ that sends $t \mapsto 1$, the {\it value semigroup} of the cusp is ${\rm S}:=v_t(f(R))$, and it is a fundamental topological invariant attached to the singularity.

\medskip
In this note we collect, contextualize, and invite the reader to improve upon, a number of qualitative and quantitative results for linear series on cuspidal rational curves. As our methods are purely local, we suspect that many of these results generalize to (suitably general) curves of higher genus. Cuspidal rational curves have appeared before in Brill--Noether theory, most notably in \cite{EH}.\footnote{We thank Steve Kleiman for reminding us of this.} Eisenbud and Harris showed that the Brill--Noether theory of rational curves with $g$ {\it simple} cusps in general points is identical to that of a general smooth curve of genus $g$; and this foundational result has far-reaching consequences for the birational geometry of $\mc{M}_g$. Our work runs in an orthogonal direction: instead of focusing on the (global) geometry of rational curves with $g$ simple cusps, we focus on the (local) geometry of rational curves with a single cusp of delta-invariant $g$. It is also natural to ask for results interpolating between our results for unicuspidal rational curves and those of {\it loc. cit.} in this context; such results would depend on as-yet-unknown dimensional-transversality properties for intersections of the generalized Severi varieties we introduce in the following subsection. 

\subsection{Topology and rationality of parameter spaces, and gonalities}
Plane curves comprise another distinguished class of algebraic curves, and their parameter spaces have attracted substantial interest. Zariski \cite{Z} first established an upper bound
for the dimension of any given component of the Severi variety $V_{d,g}$ of plane curves
of fixed degree $d$ and arithmetic genus $g$, and showed that whenever the upper bound is
achieved, a general curve in that component is nodal. His result then played an important role in Harris’ proof \cite{Ha} of the irreducibility of $V_{d,g}$.

\medskip
One upshot of Harris' work is that $V_{d,g}$ is the closure of the subvariety corresponding to rational curves with $g$ nodes. It is natural, then, to ask whether, when $n \geq 3$, generalized Severi varieties $M^n_{d,g}$ parameterizing irreducible and linearly nondegenerate degree-$d$ rational curves in $\mb{P}^n$ with arithmetic genus $g$ are themselves irreducible, and are the closures of their subloci indexing $g$-nodal curves. In \cite{CFM, CLM}, we show that in general the answer to both questions is ``no", by producing explicit Severi-type varieties of unicuspidal rational curves of dimension larger than $g$-nodal loci. A simple but important ingredient operative in our analysis is the stratification of unicuspidal rational curves according to their {\it ramification profiles} ${\bf r}$ in the preimages of their respective cusps\footnote{Here ${\bf r}$ records the nonzero orders of vanishing in the preimage of the cusp of the sections that define the morphism $\mb{P}^1 \ra \mb{P}^n$ whose image is the unicuspidal rational curve in question.}. In \cite{CLMR}, we produce a conjectural combinatorial formula for the (co)dimension of a Severi variety $M^n_{d,g; {\rm S},{\bf r}}$ of unicuspidal rational curves of degree $d$, arithmetic genus $g$, and cuspidal type $({\rm S},{\bf r})$ in $\mb{P}^n$; and we verify that our conjecture is valid for two interesting infinite classes of examples.

\medskip
These Severi varieties are generalizations of Schubert varieties; indeed, the space $M^n_d$ of degree-$d$ morphisms $\mb{P}^1 \ra \mb{P}^n$ with irreducible and linearly nondegenerate images is parameterized by an open subset of the Grassmannian $\mb{G}=\mb{G}(n,d)$, and those maps with ramification profiles ${\bf r}$ in a (variable) point determine a Schubert subvariety $M^n_{d;{\bf r}} \sub \mb{G}$. Schubert varieties, which parameterize rational curves subject to ramification profiles in points, are rational; and it is natural to wonder to what extent this property persists for Severi varieties of unicuspidal rational curves, in which nonlinear conditions {\it beyond ramification} are imposed by the additive structure of the value semigroups of the cusps. In all examples that have been studied to date, the Severi varieties $M^n_{d,g; {\rm S},{\bf r}}$ have turned out to be {\it unirational}.

\medskip
In an orthogonal vein, the {\it gonality} of a curve $C$, i.e., the lowest degree $d$ for which a degree-$d$ morphism $C \ra \mb{P}^1$ exists, is a fundamental measure of its (ir)rationality. Every morphism $C \ra \mb{P}^1$, in turn, is specified by a meromorphic function on $C$. When $C$ is a unicuspidal rational curve with value semigroup ${\rm S}$, every such function satisfies conditions imposed by arithmetic properties of ${\rm S}$.

\section{Topology of Severi varieties of unicuspidal rational curves}

\subsection{Dimensionality and unirationality conjectures}
To compute the number of algebraically independent conditions imposed on morphisms $\mb{P}^1 \ra \mb{P}^n$ of fixed degree $d$ imposed by a cusp of type $({\rm S},{\bf r})$, we carry out a combinatorial local analysis of power series in the preimage of the cusp. More precisely, we develop a combinatorial framework that enables us to explicitly characterize those conditions on power series coefficients imposed by the pair $({\rm S},{\bf r})$. The local algebra of the cusp given by the image of \eqref{ring_map} is determined by its truncation modulo $t^{2g}$, where $g$ is the (local contribution to the arithmetic) genus of the cusp. This algebraic phenomenon manifests combinatorially as the fact that any numerical semigroup ${\rm S}$ of genus $g:=\mb{N} \setminus {\rm S}$ is determined by its truncation ${\rm S} \cap [2g]$.\footnote{Crucially, no element of $\mb{N} \setminus {\rm S}$ is strictly greater than $2g$.} It quickly becomes clear that combinatorics provides a natural organization scheme for local algebra.

\medskip
Explicitly, to count conditions imposed on coefficients of morphisms $f:\mb{P}^1 \ra \mb{P}^n$ by cusps of type $({\rm S},{\bf r})$, we begin by assigning a lattice path to ${\rm S}$ from $(0,0)$ to $(g,g)$ in $\mb{Z}^2$, with horizontal and vertical segments of unit length, as follows: for each successive integer $i \in [2g]$, we proceed to the right (resp., upward) if $i \in {\rm S}$ (resp., $i \notin {\rm S}$). Any morphism $f:\mb{P}^1 \ra \mb{P}^n$ whose image has a cusp of type $({\rm S},{\bf r})$ is given locally near the cusp by an $n$-tuple of power series $f_i$ with valuations $r_i, i=1,\dots,n$, where ${\bf r}=(r_1,\dots,r_n)$; so the integers $r_i$ belong to ${\rm S}$, and index columns of the lattice square $\{(i,j): 1 \leq i,j \leq g\}$. The column corresponding to $f_i$ now contributes two types of conditions: ramification conditions, imposed by the (linear) vanishing of $j$-th order derivatives $D^j_t$, $j=0,\dots,r_i-1$ in the preimage of the cusp; and those (nonlinear) conditions beyond ramification that depend on the additive structure of ${\rm S}$.

\medskip
Unpacking conditions beyond ramification requires a closer examination of ${\rm S}$ from an arithmetic point of view. Not surprisingly, the (unique) set of minimal generators of ${\rm S}$ plays a key role. More obscure, but equally important for our purposes, are the {\it Betti} elements of the pair $({\rm S},{\bf r})$, which determine the structure of factorizations of elements in ${\rm S}$ as nonnegative linear combinations of minimal generators and the ramification orders $r_i$, $i=1,\dots,n$ \footnote{The usage of ``Betti" here is consistent with its standard usage in describing the free resolutions of ideals. Indeed, Betti elements of a finite set ${\rm T}$ of positive integers select for minimal generators of the toric ideal of the image of the rational curve determined by monomials $t^e$, $e \in {\rm T}$.}. Formally, given any finite set ${\rm T}$ of (strictly) positive integers, let $\rm{F}_{\rm T}$ (resp., $\rm{S}_{\rm T}$) denote the free monoid (resp., numerical semigroup) they generate. There is a natural projection $\pi: \rm{F}_{\rm T} \ra \rm{S}_{\rm T}$, whose fibers describe how elements of ${\rm S}_{\rm T}$ factor as nonnegative linear combinations of elements in ${\rm T}$. Decreeing that $v \sim w$ whenever elements $v, w \in \pi^{-1}(s)$ fit into a chain $v=v_1, \dots, v_n=w$ for which $\langle v_i, v_{i+1} \rangle \neq 0$ for every $i$ yields an equivalence relation on $\pi^{-1}(s)$; and $s$ is a {\it Betti} element whenever $\pi^{-1}(s)$ splits into at least two $\sim$-equivalence classes.

\medskip
The Betti elements of $({\rm S},{\bf r})$ index precisely those columns in the Dyck diagram corresponding to parameterizing functions $f_i$ that impose conditions beyond ramification. Conditions beyond ramification, in turn, arise when lowest-order terms of polynomials in $f_i$, $i=1,\dots,n$ (whose valuations belong to ${\rm S}$) are forced to vanish. In \cite{CLMR}, we formulate an explicit combinatorial conjecture to account for these conditions; while these are nonlinear in general, all are indexed by linear algebraic quantities, namely 
jumps in the ranks of matrices  derived from the Dyck diagram. Explicitly, let ${\bf r}^*=\{r_{n+1},\dots,r_{\ell}\}$ denote the set of minimal generators of ${\rm S}$ less than the conductor $c$ that do not belong to ${\bf r}$; and let $B$ denote the set of Betti elements of ${\bf r} \sqcup {\bf r^*}$ strictly less than $c$. Given $b\in B$, let $\{v_1,\ldots,v_{n_b}\}$ be any full set of representatives of $\sim$-equivalence classes of $b$; let $M_b$ denote the $(n_b-1)\times \ell$ matrix
$\begin{pmatrix}
v_2-v_1\\
\vdots\\
v_{n_b}-v_1
\end{pmatrix}
$.
Assuming that $b_1<\dots<b_p$ comprise $B$, let
$
A_i :=
\begin{pmatrix}
M_{b_1}\\
\vdots\\
M_{b_i}
\end{pmatrix}
$
and given $b=b_i\in B$, let
$\phi(b):={\rm rank}(A_i)-{\rm rank}(A_{i-1})$. Finally, let
$
B^{\pr}:=\{ b\in B\,|\,\phi(b)\geq 1\}
$ and let
${\bf r}^{\bullet}:=\{ r_i \in {\bf r}^{*} \ |\  B^{\pr}\cap (r_{i-1},r_i)\neq\emptyset\text{ or } \#\{ b\in B^{\pr}\ | \ b < r_i \} >  (i-n)-1\}$.

\begin{conj}\label{reformulated_conjecture}
Given a vector ${\bf r}=(r_1,\ldots,r_n)\in \mathbb{N}_{>0}^n$, let $\mathcal{V}_{\bf r} \subset M^n_{d,g; {\rm S}}$ be the subvariety of maps $f: \mb{P}^1 \rightarrow \mb{P}^n$ with a unique singularity that is cuspidal with semigroup ${\rm S}$ and ramification profile ${\bf r}$. Assume $d=\deg(f) \geq \max(n,2g-2)$. Then
\begin{equation}\label{general_cod_estimateBis}
{\rm cod}(\mc{V}_{\bf r},M^n_d) =\sum_{i=1}^n(r_i-i)+\sum_{s\in B} \phi(s)\rho(s)-\sum_{s\in{\bf r}^{\bullet}} \rho(s)- 1
\end{equation}
where $\rho(s)$ denotes the number of elements of $\mb{N} \setminus {\rm S}$ strictly greater than $s \in {\rm S}$.
\end{conj}

On the right hand side of \eqref{general_cod_estimateBis}, conditions beyond ramification account for \\$\sum_{s\in B} \phi(s)\rho(s)-\sum_{s\in{\bf r}^{\bullet}} \rho(s)$.
Graphically speaking, for every Betti element $s \in B$, $\rho(s)$ computes the number of squares above the Dyck path in the column corresponding to $s$, while $\phi(s)$ gives a ``multiplicity" with which that column contributes. Conjecture~\ref{reformulated_conjecture} reflects the expectation that a complete set of algebraically independent conditions beyond ramification imposed by cusps of type $({\rm S},{\bf r})$ is produced by an explicit inductive algorithm detailed in \cite{CLM}.

\begin{ex}\label{hyperelliptic_example} Let $n=4, g=7$, ${\rm S}=\langle 2,15 \rangle$ and ${\bf r}=(2,4,6,8)$;
see Figure~\ref{hyperelliptic_graphic}.
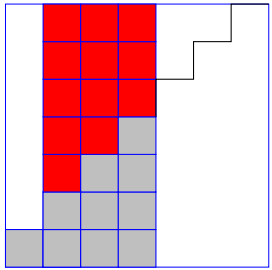
\begin{figure}\label{hyperelliptic_graphic}
\begin{tikzpicture}[scale=0.50,squarednode/.style={rectangle, draw=red!60, fill=red!5, very thick}]
\draw[blue, very thin] (0,0) rectangle (7,7);
\filldraw[draw=blue, fill=lightgray] (0,0) rectangle (1,1);

\filldraw[draw=blue, fill=lightgray] (1,0) rectangle (2,1);
\filldraw[draw=blue, fill=lightgray] (1,1) rectangle (2,2);
\filldraw[draw=blue, fill=lightgray] (1,0) rectangle (2,1);
\filldraw[draw=blue, fill=lightgray] (1,1) rectangle (2,2);
\filldraw[draw=blue, fill=red] (1,2) rectangle (2,3);
\filldraw[draw=blue, fill=red] (1,3) rectangle (2,4);
\filldraw[draw=blue, fill=red] (1,4) rectangle (2,5);
\filldraw[draw=blue, fill=red] (1,5) rectangle (2,6);
\filldraw[draw=blue, fill=red] (1,6) rectangle (2,7);

\filldraw[draw=blue, fill=lightgray] (2,0) rectangle (3,1);
\filldraw[draw=blue, fill=lightgray] (2,1) rectangle (3,2);
\filldraw[draw=blue, fill=lightgray] (2,2) rectangle (3,3);
\filldraw[draw=blue, fill=red] (2,3) rectangle (3,4);
\filldraw[draw=blue, fill=red] (2,4) rectangle (3,5);
\filldraw[draw=blue, fill=red] (2,5) rectangle (3,6);
\filldraw[draw=blue, fill=red] (2,6) rectangle (3,7);

\filldraw[draw=blue, fill=lightgray] (3,0) rectangle (4,1);
\filldraw[draw=blue, fill=lightgray] (3,1) rectangle (4,2);
\filldraw[draw=blue, fill=lightgray] (3,2) rectangle (4,3);
\filldraw[draw=blue, fill=lightgray] (3,3) rectangle (4,4);
\filldraw[draw=blue, fill=red] (3,4) rectangle (4,5);
\filldraw[draw=blue, fill=red] (3,5) rectangle (4,6);
\filldraw[draw=blue, fill=red] (3,6) rectangle (4,7);

\draw (4,4) -- (4,5) -- (5,5) -- (5,6) -- (6,6) -- (6,7) -- (7,7);

\end{tikzpicture}

\caption{Ramification conditions are in grey, while conditions beyond ramification are in red; conditions contributed by $f_i$ are in the $i$th column.}
\end{figure}
Renormalizing, write $f_i= \sum_{j=0}^{\infty} a_{i,j} t^j$ with $a_{i,j}=0$ for all $j<2i$ and $a_{i,2i}=1$; and set $F_i:= f_i-f_1^i$. A typical condition beyond ramification is that $\text{lc}(F_i)=a_{i,2i+1}-i a_{1,3}=0$, $i=2,\dots,n$, where ``lc" denotes the coefficient of the term of $F_i$ with valuation $v_t(f_i)+1=v_t(f_1^i)+1=2i+1$. We have $\text{lc}(F_i)=0$  because $2i+1 \notin {\rm S}$; and this condition is encoded by the lowest red square in column $i$. Inductively, we now ``walk" up column $i$, inducing a single new independent condition at every step. The condition encoded by the second-lowest red square is imposed by $F_i^{\ast}:= F_i- [t^{2i+2}]F_i \cdot f_1^i$. To continue walking up the column, we replace $F_j$ by $F_j^{\ast}$ and perturb by (a multiple of) a power of $f_1$; and
iterate this procedure until all elements of $\mb{N} \setminus {\rm S}$ have been exhausted.
\end{ex}

A key feature of Example~\ref{hyperelliptic_example} is that for every column $i \geq 2$ in the Dyck diagram, the expansions of the leading coefficients of $F_i$ and its inductively-created counterparts $F_i^*$ whose vanishing is imposed by $({\rm S},{\bf r})$ always contain (linear) instances of variable coefficients $a_{i,j}=[t^j]f_i$ of $f_i$ that are not present in earlier iterations of the algorithm. This implies that the corresponding Severi variety $\mc{V}_{\bf r}$ is unirational.

\begin{conj}\label{unirationality_conj}
With hypotheses and notation as in Conjecture~\ref{reformulated_conjecture}, the Severi variety $\mc{V}_{\bf r}$ of unicuspidal rational curves of type $({\rm S},{\bf r})$ is unirational.
\end{conj}
 Conjecture~\ref{reformulated_conjecture} also naturally begs the following two questions:
\begin{enumerate}
\item[i.] Are there explicit upper and lower polynomial bounds in $n$ and $g$ for the right hand side of \eqref{general_cod_estimateBis}?
\item[ii.] What are quantitatively precise versions of Conjecture~\ref{reformulated_conjecture} for distinguished infinite families of types $({\rm S},{\bf r})$?
\end{enumerate}
The following two sections aim to shed light on each of these. 

\subsection{Lattice point counts and numerology}
Recall that a point $P$ of a curve $C$ is \emph{Gorenstein} whenever $\ww_P$ is a free $\oo_P$-module; and $C$ is \emph{Gorenstein} when this holds for every $P$. In this case, $\ww$ is an invertible sheaf, and the value semigroup ${\rm S}$ at every $P$ is {\it symmetric}, i.e., for every $x \in {\mb N}$, either $x$ or $c-1-x$ belongs to ${\rm S}$.\footnote{In fact, $P \in C$ is Gorenstein if and only if its value semigroup ${\rm S}$ is symmetric.} Whenever ${\rm S}$ is symmetric, in turn, we may rewrite the right hand side of \eqref{general_cod_estimateBis} in terms of counts of interior lattice points of polytopes in $\mb{R}^e$, where $e=e({\rm S})$ is the number of minimal generators of ${\rm S}$. Hereafter we will focus on two infinite families of types $({\rm S},{\bf r})$ associated with very differently-structured Betti sets $B=B({\rm S},{\bf r})$. The first is comprised of semigroups ${\rm S}$ that are {\it hyperelliptic}, i.e., for which $2 \in {\rm S}$, and ramification profiles ${\bf r}=(2,4,\dots, 2n)$; while the second is comprised of {\it supersymmetric} semigroups ${\rm S}$ minimally generated by products $\frac{a_1 \cdots a_n}{a_i}$, $i=1,\dots,n$ and ramification profiles ${\bf r}$ that are minimal generating sets.\footnote{Here we assume $a_1<\dots<a_n$ are pairwise relatively prime positive integers.} These two families behave very differently with regard to their Betti elements. Indeed, in the first case, the $(n-1)$ even numbers between 4 and $2n$ belong to $B({\rm S},{\bf r})$. On the other hand, supersymmetric semigroups are precisely those numerical
semigroups whose sets of Betti elements are singletons \cite[Thm 12]{GOR}; so cases in which ${\rm S}$ is
supersymmetric and ${\bf r}$ comprises a set of minimal generators are minimal with
respect to $\#B({\rm S}, {\bf r})$.

\begin{thm}[\cite{CLMR}, Thm 2.1 and Thm 4.2]\label{dimensionality_theorem}
Let $d$, $g$, and $n$ be positive integers for which $n \leq 2g$ and $d \geq \max(n, 2g-2)$. 
\begin{itemize}
\item The subvariety of $M^n_d$ parameterizing unicuspidal rational curves with ${\rm S}=\langle 2,2g+1\rangle$ and ${\bf r}=(2,4,\dots,2n)$ is unirational and of codimension
$(n-1)g$.

\item Let $n=3$; and let $a_1$, $a_2$, and $a_3$ be pairwise relatively prime positive
integers. The subvariety of $M^3_d$ of unicuspidal rational curves with ${\rm S} = \langle a_1a_2, a_1a_3, a_2a_3 \rangle$ and ramification profile ${\bf r}=(a_1a_2, a_1a_3, a_2a_3)$ is unirational and of codimension $2\rho(a_1a_2a_3)+a_1a_2+a_1a_3+a_2a_3-7$.
\end{itemize}
\end{thm}

Theorem~\ref{dimensionality_theorem} confirms Conjectures~\ref{reformulated_conjecture} and \ref{unirationality_conj} for our distinguished hyperelliptic family of examples and in supersymmetric cases with embedding dimension $n=3$. Note that the supersymmetric codimension formula involves $\rho(a_1a_2a_3)$, which is difficult to identify explicitly in general. On the other hand, we may {\it approximate} $\rho(a_1a_2a_3)$ by re-interpreting it as a discrete volume of a 3-dimensional rational polytope, as follows. Every supersymmetric semigroup ${\rm S}=\langle a_1a_2,a_1a_3,a_2a_3 \rangle$ is symmetric, with conductor $c=2a_1a_2a_3-(a_1a_2+a_1a_3+a_2a_3)+1$; as a result, elements of $\mb{N} \setminus {\rm S}$ strictly greater than $a_1a_2a_3$ are in bijection with elements of ${\rm S}$ strictly less than $a_1a_2a_3-(a_1a_2+a_1a_3+a_2a_3)$. Elements of ${\rm S}$ strictly less than the Betti element $a_1a_2a_3$ factor uniquely as anonnegative linear combination of minimal generators; so $\rho(a_1a_2a_3)$ computes the number of lattice points inside the simplex $\Delta$ with vertices $(0,0,0)$, $(a_3-1-\frac{a_3}{a_1}-\frac{a_3}{a_2},0,0)$,$(0,a_2-1-\frac{1}{a_1}-\frac{1}{a_3},0)$, and $(0,0,a_1-1-\frac{1}{a_2}-\frac{1}{a_3})$. On the other hand, the discrete volume of $\Delta$ is approximated by its Euclidean (``actual") volume; a precise upper bound is given in \cite{YZ}; and this, in turn, allows us to conclude in many cases that the codimension inside $M^3_d$ of the Severi variety of supersymmetric unicuspidal rational curves is strictly less than $g$, the codimension of the $g$-nodal locus \cite{CKV}.

\begin{spec}
Let $d$, $g$, and $n$ be positive integers for which $n \leq 2g$ and $d \geq \max(n, 2g-2)$. Given a numerical semigroup ${\rm S}$ of genus $g$, the corresponding Severi variety $M^n_{d,g; {\rm S}}$ of unicuspidal rational curves with semigroup ${\rm S}$ is of codimension at most $(n-1)g$, with equality if and only if ${\rm S}=\langle 2,2g+1 \rangle$.
\end{spec}

\subsection{Certifying dimensionality and unirationality}

Certifying Conjecture~\ref{reformulated_conjecture} in the hyperelliptic and superelliptic cases is  the focus of \cite{CLMR}; that article also establishes a template for certifying and/or refuting Conjecture~\ref{reformulated_conjecture} in general. The conjectural premise here is that the Betti elements $B({\rm S},{\bf r})$ determine a complete set of algebraic relations for the local algebra of an arbitrary cusp of type $({\rm S},{\bf r})$. Whether or not this premise is valid is algorithmically checkable in any particular case via a step-by-step procedure, starting with a ``universal" element $F$ of the corresponding local algebra, and alternating between imposing
vanishing on coefficients of $F$ in orders corresponding to elements of ${\rm S}$ and deducing
polynomial vanishing conditions imposed by terms whose orders are
elements of $\mb{N} \setminus {\rm S}$. For hyperelliptic and supersymmetric families, the Betti sets $B({\rm S},{\bf r})$ are explicit and their orderings are particularly simple; it remains to generalize the arguments of \cite{CLMR} in order to accomodate arbitrary Betti structures.

\section{Canonical models}
\subsection{The canonical model of a singular curve}

Let $C$ be an integral and complete curve over an algebraically closed field of (arithmetic) genus $g$ with structure sheaf $\oo=\oo_C$, dualizing sheaf $\ww=\ww_C$, and field of meromorphic functions $k(C)$. 
Given a coherent sheaf $\mathcal{F}$ on $C$ and a morphism $\varphi:X\to C$ from a scheme $X$ to $C$, let
$\oo_{X}\mathcal{F}:=\varphi^* \mathcal{F}/\rm{Torsion}(\varphi^*\mathcal{F})$ and $\fff^n:=\rm Sym^n\fff/\rm Torsion (\rm {Sym}^n\fff)$.

\medskip
Following \cite{S}, a \emph{linear series of rank $r$ and degree $d$ on $C$} is a set
$\sys(\fff ,V):=\{x^{-1}\aaa\ |\ x\in V\setminus 0\}$,
where $\mathcal{F}$ is a torsion-free sheaf of rank $1$ and degree $d$ on $C$ and $V \sub H^{0}(\aaa )$ is a vector subspace of dimension $r+1$. 
Here $V\subset k(C)$ and $x^{-1}\fff$ is the sheaf characterized locally by $(x^{-1}\fff)(U)=x^{-1}\fff(U)$ for any open set $U\subset C$.

\medskip
The \emph{canonical model} of $C$, introduced by Rosenlicht \cite{R}, is the image $C'$ of the curve $C$ under the morphism 
$\kao:\overline{C}\rightarrow\pp^{g-1}$ defined by the linear series $\sys(\oo_{\cb}\ww,H^0(\ww))$. Rosenlicht's main theorem establishes that whenever
$C$ is nonhyperelliptic, its normalization map factors through a morphism $C'\to C$; and it holds for curves $C$ that are arbitrarily singular, including non-Gorenstein curves. In \cite{KM}, the authors re-interpreted Rosenlicht's result as an isomorphism between the canonical model and the \emph{canonical blowup} $\widehat{C}:=\text{Proj}(\oplus\,\ww ^n)$. 

\medskip
A cornerstone of the classical study of algebraic curves is {\it Max Noether's theorem}, which establishes that for any smooth curve $C$, the canonical (multiplication) morphism $\mu_n: \text{Sym}^n\ H^0(\ww) \ra H^0(\ww^n)$ is surjective for every positive integer $n$, the proof of which relies on auxiliary surjections $\text{Sym}^n\ H^0(\ww) \ra H^0(\oo_{C'}(n))$ through which each $\mu_n$ factors. Miraculously, the analogous maps $\mu_n$ are also surjective when $C$ is singular, even though $\oo_{C'}(n)$ may not coincide with $\ww^n$ in general.\footnote{Indeed, we have $\oo_{C'}(n)=\ww^n$ for every $n \geq $ if and only if $C$ is either Gorenstein or 
{\it nearly Gorenstein} in the sense of \cite{KM}.}

\medskip
A key operative ingredient in the proof of Max Noether's theorem in the singular case 
is the 
fact that for every $P\in C$, 
we have $v(\wwp)=\kk$, where $v:R \ra {\rm S}$ denotes the natural order-of-vanishing valuation on the local ring $R$ of $C$ in $P$, and 
$\kk:=\{ a\in\zz\ |\ c-1 -a\not\in\sss\}$; and that the multiplication morphism $\mu_n$ admits a valuation-theoretic analogue involving $\kk$.
In \cite{Mt}, 
surjectivity of (all of) the multiplication maps $\mu_n$ in the cuspidal case was ultimately reduced to (and deduced from) the statement that for every numerical semigroup $\sss$ with conductor $c$, every $n\in \{c,c+1,\ldots,2c-3\}$ may be written as a sum $n=k_1+k_2$ with $k_i\in\kk$ and $k_i<c$, $i=1,2$. In \cite{CFM}, this numerical statement was generalized for semigroups of rank greater than one (arising from multibranch singularities) and proved for semigroups of rank two (from singularities with two branches); and in \cite{GM}, the generalized 
statement was proved for semigroups of arbitrary rank, yielding Max Noether's 
theorem for arbitrary integral curves.

\subsection{Gonalities from canonical models} 

The {\it gonality} of a smooth curve $C$ is 
the smallest $d$ for which there exists a degree-$d$ cover $C\to\pum$. If $C$ is non-hyperelliptic, it is intimately related to the projective geometry of 
the canonical embedding of $C$. Indeed, 
$C$ is $d$-gonal if and only if the canonical image of $C$ lies on a $(d-1)$-fold scroll. 
The expectation is that this phenomenon persists for singular curves, provided we replace the canonical image by the canonical model and define gonality to be the smallest
$d$ for which $C$ admits a 
rank one torsion-free sheaf $\aaa$ with $\deg(\aaa)=d$ and $h^0(\aaa)\geq 2$.\footnote{In the terminology of \cite{RSt}, this means that $C$ admits a $g_d^1$ with a {\it non-removable} base point.}
In \cite{NMM,LMS} this was established for rational monomial curves via a combinatorial argument that we describe next.

\medskip
We begin by identifying the canonical model $C'$ in terms of $C$. So assume $C$ is given explicitly by  
\begin{equation}
\label{equrmc}
C=(t^{a_0}:t^{a_1}:\cdots:t^{a_{n-1}}:t^{a_n})\subset \mathbb{P}^n
\end{equation}
with $0=a_0<a_1<\ldots < a_n$. $C$ then has (at most) two singularities supported in $P_1=(1:0:\cdots:0)$ and 
$P_2=(0:0:\cdots : 1)$ with semigroups $\sss_{P_1}=\langle a_i\rangle_{i=1}^{n}$ and $\sss_{P_2}=\langle a_n-a_i\rangle_{i=0}^{n-1}$, respectively. Let $G_{P_i}:= \mb{N} \setminus \sss_{P_i}$, 
and set $\delta_i=\#(G_i)$; the arithmetic genus of $C$ is $g=\delta_1+\delta_2$. Let $\gamma$ be the {\it Frobenius number} of $\sss_{P_1}$, i.e., the largest integer in $G_1$. The authors of \cite{LMS} establish that
\begin{equation}
\label{equcmo}
C'=(1:t^{b_2}:\cdots :t^{b_{\delta_P}}:t^{c_1}:\cdots:t^{c_{\delta_Q}})\subset\mathbb{P}^{g-1}
\end{equation}
where $\{0,b_2,\dots,b_{\delta_P}\}=\gamma-G_{P_1}$ and $\{c_1,\dots,c_{\delta_Q}\}=\gamma+G_{P_2}$. The proof relies heavily on the equality $v(\wwp)=\kk$ applied to the singularities $P=P_i$, $i=1,2$.

\medskip
We also need 
to identify when a monomial curve lies on a scroll. The relevant condition, proved in \cite{NMM}, is that $(1:t^{a_1}:\ldots :t^{a_{N}})\subset \pp^{N}$ lies on a $d$-fold scroll if and only if 
the set $\{0=a_0,a_1,\ldots,a_{N}\}$ may be partitioned into $d$ subsets whose elements are in arithmetic progression with the same common difference $\mu$.

\medskip

Finally note that any sheaf computing the gonality $d$ of $C$ is of the form $\oo_C\langle 1,t^{\mu}\rangle$ for some $\mu \in\mathbb{N}$. In \cite{LMS} the authors prove that the exponents $b_i$ and $c_i$ of equation \eqref{equcmo} may be partitioned into $d-1$ arithmetic progressions each with common difference $\mu$, 
and that the resulting $(d-1)$-fold scroll on which $C'$ lies is of minimal dimension.

\begin{ex}
\emph{
Assume $C \sub \mb{P}^n$ is a unicuspidal monomial rational curve. Let $P_1=(1:0:\cdots:0)$ and $P_2=(0:\cdots:0:1)$ as before, assume the cusp of $C$ is supported in $P_1$, and let $\sss :=\sss_{P_1}$. 
For every $\mu\in\nn$, we have
\begin{equation}
\label{equssm}
\deg\oo_C\langle 1,t^{\mu}\rangle = \#(((\sss\cup(\sss+\mu))\setminus \sss) + \mu.
\end{equation}
Indeed, 
$\aaa:=\oo_C\langle 1,t^{\mu}\rangle$ is of degree zero away from $P_1$ and $P_2$.  On the other hand, $t^{-1}$ is a local parameter at $P_2$, so $\deg_{P_2}(\aaa)=\mu$; while the first summand on the right-hand side of \eqref{equssm} is $\deg_{P_1}(\mc{F})=\deg_{P_1}(\mc{O}_{P_1} +t^{\mu}\mc{O}_{P_1})$.}

\medskip
\emph{Now suppose $C=(1:t^5:t^7:t^8)\subset \mathbb{P}^3$. 
In this case, $\sss=\langle 5,7,8 \rangle$ is of genus 7. 
We will show that $C$ is 4-gonal, and that its gonality is computed by the sheaf $\aaa=\oo_C\langle 1,t^2\rangle$. Indeed, that $\deg(\aaa)=4$ follows from \eqref{equssm} and the fact that $(((\sss\cup(\sss+2))\setminus \sss)=\{2,9\}$. On the other hand, the gonality of a unicuspidal rational monomial curve 
is always realized by a sheaf of the form $\mc{O}\langle 1,t^{\mu}\rangle$ \cite{FM}. 
Applying \eqref{equssm} with $\mu =1$, 3 and 4 yields degrees $5,5$ and $7$ respectively; meanwhile, for every $\mu\geq 5$, it follows immediately from the right-hand side of \eqref{equssm} that $\deg\oo_C\langle 1,t^{\mu}\rangle\geq \mu$.} 

\medskip
\noindent \emph{Here the canonical model of $C$ is
$$
C'=(1:t^2:t^5:t^7:t^8:t^9:t^{10}) 
\subset \mathbb{P}^6.
$$
Partitioning powers of $t$ into sets $\{0,2\}$, $\{5,7,9\}$, and $\{8,10\}$ each with common difference $2$, we see that $C'$ lies on the 3-fold scroll 
\[
S_{112}=
\bigg(
\begin{array}{cccc}
x_0 & x_4 &  x_2 & x_3 \\
x_1 & x_6 & x_3 & x_5
\end{array}
\bigg)
\]
in which the subscript indicates that $S=\mb{P}(\oo_{\mb{P}^1}(1)\oplus\oo_{\mb{P}^1}(1)\oplus\oo_{\mb{P}^1}(2))$ as a $\mb{P}^1$-bundle. Moreover, $S_{112}$ is of minimal dimension among scrolls containing $C'$ because $\{0,2,5,7,8,9\}$ cannot be partitioned into two subsets comprising arithmetic progressions with the same common difference.
}
 \end{ex}

\subsection{A conjectural upper bound for covering gonality}

It is natural to ask for an upper bound on the gonality of a unicuspidal rational curve $C$; this is the focus of work in progress of the authors together with V. Lara and N. Galdino. Here we focus on a type of gonality distinct from that of the previous subsection, the {\it covering gonality}, which we denote by $\gon_F(C)$: namely, the smallest $d$ for which a $d$-fold cover $C\to\pum$ exists.\footnote{Equivalently, $\gon_F(C)$ is the smallest $d$ for which $C$ admits a basepoint-free $g^1_d$.}

\begin{conj} A Gorenstein unicuspidal rational curve $C$ has covering gonality $\gon_F(C)\leq \lceil 
\frac{g+1}{2} \rceil$.
\end{conj}\label{covering_gonality_conj}

We conclude this note by giving a heuristic argument in favor of Conjecture~\ref{covering_gonality_conj}. To this end, note that $C$ admits a degree $k$ morphism to $\pum$ if and only if 
there is some point $P \in C$ of multiplicity $m=m(P)$ and a meromorphic function $f\in\op$ of the form
$$
f:=\frac{d_mt^m+d_{m+1}t^{m+1}+\ldots +d_kt^k}{1+r_1t+\ldots +r_kt^k}
$$
with either $d_k\neq 0$ or $r_k\neq 0$.

\medskip \noindent
Now suppose that
$$
C=(F_0 : F_{1} : \ldots : F_{g-1}) \subset \mathbb{P}^{g-1}
$$
with $F_i\in k[t]$ and $\deg(F_i)\leq 2g-2$; and suppose, furthermore, that $v_t(F_i)$ is the $(i+1)$-st element of the value semigroup of the singularity of $C$ in $P$. We need to find $d_i$, $r_i$ and $c_i$ for which
$$
\frac{d_mt^m+d_{m+1}t^{m+1}+\ldots +d_kt^k}{1+r_1t+\ldots +r_kt^k}=\frac{c_{1}F_{1}+\cdots+c_{g-1}F_{g-1}+t^{2g}u(t)}{F_0}.
$$
This amounts to solving a nonlinear system of equations in which there are
\[
n_v=\underbrace{k-m+1}_{d_i}+\underbrace{k}_{r_i}+\underbrace{g-1}_{c_i}
\]
variables and
$
n_e=2g-m
$
equations, with each $t^i$, $i=m,\dots,2g-1$ contributing an equation. Naively, we'd expect the system to be solvable whenever $n_v \geq n_e$. However, when $g=4$ and $m=k=2$ we have $n_v = n_e$, and it may be checked that $k$-gonal curves 
comprise a proper closed subvariety of $M^6_{4,\langle 2,9\rangle,(2,4,6)}$. 
More generally, we anticipate our nonlinear system of equations to be solvable (for {\it all} curves in the relevant Severi variety) whenever $n_v>n_e$, i.e., whenever $k\geq \frac{g+1}{2}$. 

\end{document}